\documentclass[reqno,11pt]{amsart}

\newtheorem{theorem}{Theorem}[section]
\newtheorem{lemma}[theorem]{Lemma}
\newtheorem{corollary}[theorem]{Corollary}

\theoremstyle{definition}
\newtheorem{assumption}[theorem]{Assumption}
\newtheorem{definition}[theorem]{Definition}
 
\theoremstyle{remark}
\newtheorem{remark}[theorem]{Remark}

 \makeatletter 
 \def\dashint{%
 \operatorname%
 {\,\,\text{\bf--}\kern-.98em\DOTSI\intop\ilimits@\!\!}}
\def\dashnorm{\,\,\text{\bf--}\kern-.5em\|}
\def\ninf{\qopname\relax\@empty{inf\phantom{p}\!\!\!}}

 \makeatother

\newcommand\bR{\mathbb{R}}

\newcommand\cB{\mathcal{B}}
\newcommand\cF{\mathcal{F}}

\newcommand{\osc}{{\rm osc}\,}
\newcommand{\loc}{{\rm loc}\,}

 \newcommand{\mysection}[1]{\section{#1}
 \setcounter{equation}{0}}
 
\newcommand{\nlimsup}{\operatornamewithlimits{\overline{lim}}}
\newcommand{\nliminf}{\operatornamewithlimits{\underline{lim}}}

\begin{document}

\title[Strong solutions of It\^o's equations]
{On strong solutions of It\^o's equations
with $D\sigma $ and
   $b $ in Morrey classes containing $L_{d}$}
\author{N.V. Krylov}

\email{nkrylov@umn.edu}
\address{127 Vincent Hall, University of Minnesota, Minneapolis, MN, 55455}
 
 

\begin{abstract} 
We consider It\^o uniformly nondegenerate equations
with time independent coefficients, the diffusion
coefficient in $W^{1}_{2+\varepsilon,\loc}$, and the drift in a Morrey class containing
$L_{d}$. We prove the unique strong solvability
in the class of admissible solutions
for any starting point. The result is new even
if the diffusion is constant.
\end{abstract}

\maketitle

\mysection{Introduction}
                                                  \label{section 3.11.1}
Let $\bR^{d}$ be a $d-$dimensional Euclidean space of points
$x=(x^{1},...,x^{d})$ with $d\geq3$. Let $(\Omega,\cF,P)$ be a 
complete probability space,
let $\{\cF_{t}\}$ be an increasing filtration of 
$\sigma$-fields $\cF_{t}\subset \cF$, that are complete.
Let  
$w_{t}$ be a $d_{1}$-dimensional Wiener process relative to
$\{\cF_{t}\}$, where $d_{1}\geq d$.

Assume that on $\bR^{d}$ we are given $\bR^{d}$-valued Borel
functions $b,\sigma^{k}=(\sigma^{ik})$, $k=1,...,d_{1}$. We are going
to   investigate
the equation
\begin{equation}
                                         \label{6.15.2}
 x_{t}=x +\int_{0}^{t}\sigma^{k}(x_{s})\,dw^{k}_{s}+
\int_{0}^{t}b(x_{s})\,ds,
\end{equation}
where and everywhere below the summation 
over repeated indices is understood.

We are interested in the so-called strong solutions, that is
solutions such that, for each $t\geq0$, $x_{t}$ is $\cF^{w}_{t}$-measurable,
where $\cF^{w}_{t}$ is the completion of $\sigma(w_{s}:s\leq t)$. We present
sufficient conditions for the equation to have a strong solution
and also for the solution
to be unique (strong uniqueness) generalizing those in  article \cite{KrStr} which
should be read first in order to be able to follow  what we are doing here. For this
reason we restrict ourselves to a quite short
introduction without repeating the one from
\cite{KrStr}.

A very rich literature
on the {\em weak\/} uniqueness problem for \eqref{6.15.2}     is beyond
  the scope of this article.  
Concerning the strong solutions
in the {\em time dependent\/} case probably the best results belong to R\"ockner and Zhao \cite{RZ_21}, where, among very many other things, they prove existence and uniqueness of strong solutions of
equations like \eqref{6.15.2} with $b=b(t,x)$
and constant $\sigma^{k}$. 
We refer the reader to \cite{RZ_21} also for
a very good review  of the history of the problem.

The coefficient  $b(t,x)$ in \cite{RZ_21}
belong to the space 
$L_{q}([0,\infty),L_{p}(\bR^{d}))$
with  
$$
p,q\in[2,\infty],\quad \frac{d}{p}+\frac{2}{q}=1.
$$
In case $q=\infty$ they additionally assume that $b(t,\cdot)$ is a continuous $L_{d}$-valued function. In our case $b$ is independent
of $t$ and still our results are not covered by 
\cite{RZ_21} not only because we have
variable $\sigma^{k}$'s but also because our
$b$ is generally not in $L_{d}$. 

Here is an example in which we prove existence of strong solutions. Take $d=3,d_{1}=12$, and
for some numbers $\alpha,\beta,\gamma\geq0$
let $\sigma^{k}$ be the $k$th column of the
matrix given by  ($0/0:=3^{-1/2}$)
$$
\begin{pmatrix}\alpha & 0 & 0\\
0 & \alpha & 0\\
0 & 0 & \alpha
\end{pmatrix}+
\frac{\beta}{|x|}
\begin{pmatrix}
x^{1} &  x^{2} & x^{3} & 0 &  0   & 0   & 0  & 0  & 0 \\
0 & 0 & 0 & x^{1}  & x^{2}   & x^{3}   & 0  & 0  & 0 \\
0 & 0 & 0 & 0  & 0   & 0  & x^{1}& x^{2} & x^{3} 
\end{pmatrix} ,
$$
$$
b(x)=-\frac{\gamma}{|x|}\,\frac{x}{|x|}I_{0<|x|\leq 1}+\hat b(x),
$$
where $\hat b$ is, for instance, bounded.
Our result shows that if $\alpha=1$ and
$\beta$ and $\gamma$ are sufficiently small,
then \eqref{6.15.2} has a strong solution.
By the way, if $\alpha=\gamma=0$ and $\beta=1$, there exist
strong solutions of \eqref{6.15.2}
only if the starting point $x\ne0$ (see \cite{KrStr}). In case $\alpha=1$ and $\beta=0$
strong solutions exist only if $\gamma$ is sufficiently small. Observe that for $\beta\ne 0$
and $\gamma\ne0$ we have $D\sigma^{k},b\in L_{d-\varepsilon,\loc}$ for any $\varepsilon\in(0,1)$
but not for $\varepsilon=0$. Recall  that the case
that $D\sigma^{k},b\in L_{d,\loc}$ is investigated
in \cite{KrStr}, the results and methods of which
are referred to very many times in this article.

Above we were talking only about the {\em existence\/} of strong solutions. The issue of their uniqueness
is more delicate and we were able to prove
uniqueness only of ``admissible" strong solutions
(which are shown to exist).

We conclude the introduction by some notation.
We set $u_{x}=Du$ to be the gradient of $u$, $u_{xx}$ to be
the matrix of its second-order derivatives,
$$
D_{x^{i}}u=D_{i}u=u_{x^{i}}=\frac{\partial}{\partial x^{i}}u,\quad
u_{x^{i}\eta^{j}}=D_{x^{i}\eta^{j}}u=
D_{x^{i}}D_{\eta^{j}}u,
$$
$$
 \partial_{t}u=\frac{\partial}{\partial t}u,\quad
u_{(\eta)}=\eta^{i}u_{x^{i}}.
$$
If $\sigma(x)=(\sigma^{i}(x))$ is vector-valued
(column-vector), by $D\sigma=\sigma_{x}$ we mean the matrix
whose $ij$th element is $ \sigma_{x^{j}}^{i}$. If $c$ is a matrix
(in particular, vector), we set $|c|^{2}=\text{tr}\,cc^{*}$
($=\text{tr}\,c\bar c^{*}$ if $c$ is complex-valued).

For $p\in[1,\infty)$ by $L_{p}$ ($L_{p}(\Gamma)$) we mean the space 
of Borel (perhaps complex- vector- or matrix-valued)
 functions on $\bR^{d}$ (on $\Gamma\subset \bR^{d}$) with finite norm given by
$$
\|f\|_{L_{p}}^{p}=\int_{\bR^{d}}|f(x)|^{p}\,dx
\quad \Big(\|f\|_{L_{p}(\Gamma)}^{p}=\int_{\Gamma}|f(x)|^{p}\,dx\Big).
$$
By $W^{2}_{p}$ we mean the space 
of Borel functions $u$ on $\bR^{d}$ whose Sobolev derivatives 
$u_{x}$ and $u_{xx}$ exist and $u,u_{x},u_{xx}\in L_{p}$.
The  norm in $W^{2}_{p}$ is given by
$$
\|u\|_{W^{2}_{p}}=\|u_{xx}\|_{L_{p}}+\|u \|_{L_{p}}.
$$ 
Similarly  $W^{1}_{p}$ is defined. As usual, we write $f\in L_{p,\loc}$
if $f\zeta\in L_{p}$ for any $\zeta\in C^{\infty}_{0}$ ($=C^{\infty}_{0}
(\bR^{d})$).  

If a Borel $\Gamma\subset \bR^{d}$, by $|\Gamma|$ we mean its Lebesgue
measure,
$$
\dashint_{\Gamma}f(x)\,dx=\frac{1}{|\Gamma|}
\int_{\Gamma}f(x)\,dx.
$$
  Introduce
$$
B_{R}(x)=\{y\in\bR^{d}:|x-y|<R\},\quad B_{R} =B_{R}(0)
$$
and let $\cB_{R}$ be the collection of balls of radius $R$.

In the proofs of our results
we use various (finite) constants called $N$ which
may change from one occurrence to another
and depend on the data only in the same way
as it is indicated in the statements
of the results.
 
\mysection{Main results}

 Fix numbers $\delta\in(0,1)$
and $R_{0}$, $\|  b \|$,    $\|  D\sigma   \|\in(0,\infty)$.
Below   $q_{0} \in(2 ,d]$, $d\geq q> d_{0}\vee (d/2+1)$ (hence $d\geq3$), where
$d_{0}=d_{0}(d,\delta)\in(d/2,d)$
is taken from \cite{Kr_21_1}.
Set $a^{ij}=\sigma^{ik}\sigma^{jk}$, $a=(a^{ij})$.

\begin{definition}
                    \label{definition 11.14.1}
By an admissible solution of \eqref{6.15.2}
we mean any solution $x_{t}$ such that
there exists $p\in(d/2+1,q)$ and for each
$T\in(0,\infty)$ there exists a constant $N_{T}$
such that for any nonnegative Borel $f(t,x)\geq0$
\begin{equation}
                           \label{11.14.5}
E\int_{0}^{T}f(t,x_{t})\,dt\leq
N_{T}\|f\|_{L_{p}((0,T)\times\bR^{d})}.
\end{equation}
\end{definition}

\begin{assumption}
                                                  \label{assumption 3.1.1}

For any $x$ the eigenvalues of $a(x)$
lie between $\delta$ and $\delta^{-1}$,
    $b\in L_{q}$, and for any ball $B$ of radius $\rho\leq R_{0}$
\begin{equation}
                              \label{11.20.1}
 \quad\Big(\dashint_{B } |b |^{q}dx\Big)^{1/q}\leq \|  b \|  \rho^{-1}.
\end{equation}
 
\end{assumption}

\begin{assumption}
                      \label{assumption 11.3.1}
For any ball $B$ of radius $\rho\leq R_{0}$
$$
\Big(\dashint_{B } |D\sigma  |^{q_{0}}dx\Big)^{1/q_{0}}\leq \|  D\sigma  \|  \rho^{-1},\quad |D\sigma|^{2}=\sum_{i,j,k}
|\sigma^{ik}_{x^{j}}|^{2}.
$$

\end{assumption}

\begin{remark}
The case $q=q_{0}=d$ considered in
\cite{KrStr} is not excluded and in this
case one can take $\|b\|$  
as small as one likes on the account of taking $R_{0}$
sufficiently small.
  This follows
from H\"older's inequality and the fact  that
$$
 \quad\Big(\dashint_{B } |b |^{d}dx\Big)^{1/d}= N(d)\Big(\int_{B } |b |^{d}dx\Big)^{1/d}\rho^{-1}.
$$
The same goes about $\|D\sigma\|$, since
$D\sigma^{k} \in L_{d}$ in \cite{KrStr}.
Adding to this that under the conditions
in \cite{KrStr} all solutions are admissible
(see Section 4 in \cite{KrStr}) we see that the main result of \cite{KrStr} about the existence
and uniqueness of strong solutions follows from the results of the present article.

 It is also worth noting that, generally,
condition \eqref{11.20.1} with $d>q\geq d-1$ does not imply that
$b\in L_{q+ \varepsilon,\loc}$, no matter how small $\varepsilon>0$ is. Here is
an example. Take $r_{n}>0$, $n=1,2,...$, such that
the  sum of $\rho_{n}:=r_{n}^{d-q}$ is $1/2$, let $e_{1}$ be the first
basis vector, and set $b(x)=|x|^{-1}
I_{|x|<1}$, $x_{0}=1$,
$$
x_{n}=1-  2\sum_{1}^{n}r_{i}^{d-q},\quad 
c_{n}=(1/2)(x_{n}+x_{n-1})
$$
$$
  b_{n}(x)=r_{n}^{-1}b\big(r_{n}^{-1}
(x-c_{n}e_{1})\big),\quad b=\sum _{1}^{\infty}b_{n}.
$$
Since $r_{n}\leq 1$ and $d-q\leq 1$,  the supports of $b_{n}$'s are disjoint and
for $p>0$
$$
\int_{B_{1}}b^{p}\,dx=\sum _{1}^{\infty}\int_{\bR^{d}}b_{n}^{p}\,dx=N(d,p)\sum_{1}^{\infty}r_{n}^{d-p}.
$$
According to this we take the $r_{n}$'s so that
the last sum diverges for any $p>q$.
Then observe that for any $n\geq 1$ and any ball $B$
of radius $\rho$
$$
 \int_{B }  b_{n} ^{q}dx \leq N(d) \rho^{d-q} .
$$
Also, if the intersection of $B$ with $\bigcup B_{r_{n}}(c_{n})$
is nonempty, the intersection
 consists of some $B_{r_{i}}(c_{i})$, $i=i_{0},...,i_{1}$, and $B\cap B_{r_{i_{0}-1}}(c_{i_{0}-1})$ if $i_{0}\ne 0  $ and 
$B\cap B_{r_{i_{1}+1}}(c_{i_{1}+1})$. 
In this situation
$$
 \sum_{i=i_{0} }^{i_{1} }
\rho_{i} \leq 2\rho,
$$
and therefore,
$$
 \int_{B }  b  ^{q}\,dx=N(d)\sum_{i=i_{0}}^{i_{1}}r_{i}^{d-q}+\int_{B }  b_{i_{0}-1}  ^{q}dx
+\int_{B }  b_{i_{1}+1}  ^{q}dx \leq N(d)(\rho
+\rho^{d-q}),
$$
where the last term is less than $N(d)\rho^{d-q}$
for $\rho\leq 1$ and this yields  just a different form of \eqref{11.20.1}.

\end{remark}
\begin{theorem}
                        \label{theorem 11.15.1}
Under the above assumptions
there is a constant $N_{0}
\geq1$,
depending only on   $d$,   $\delta$,    $q$, $q_{0}$,  
such that if  
\begin{equation}
                                                        \label{11.15.1}
N_{0}(\|D\sigma\|+\|b\|)\leq 1,
\end{equation}
then \eqref{6.15.2} has a unique
admissible strong solution, and each admissible
solution of \eqref{6.15.2} is strong (thus coinciding with the unique one).
\end{theorem}

\begin{remark}
The solution of \eqref{6.15.2} depends on 
the starting point $x$: $x_{t}=x_{t}(x)$.
In the author's opinion, it is unlikely
that $x_{t}(x)$ is H\"older continuous
in $x$ with exponent as close to 1 as we wish
unless $q_{0}=q=d$.

How small one can take $q$ is also an interesting
question. Most likely, as $\delta\downarrow 0$,
$d_{0}\uparrow d$ and this also forces $q\uparrow d$. But if $d=d_{1}$ and we are dealing with the equation
$$
x_{t}=x+w_{t}+\int_{0}^{t}b(x_{t})\,dt
$$
our conjecture is that one can allow $q>d/2$
as close to $d/2$ as one likes and still 
have Theorem \ref{theorem 11.15.1} valid.
\end{remark}
 
\mysection{An analytic semigroup}
                                                  \label{section 3.11.2}
In this section Assumption \ref{assumption 3.1.1}
is supposed to be satisfied but
with less restrictions on $q$: namely
$q\in (1,d]$.
Assumption
\ref{assumption 11.3.1} is replaced with a   weaker Assumption
\ref{assumption 2.20.1} which comes after some preparations.

 Introduce the uniformly  elliptic operators
  $$
Lu( x)=(1/2) a^{ij}( x) u_{x^{i}x^{j}} ( x)+
b^{i}( x) u_{x^{i}}( x)\quad (a=\sigma\sigma^{*}),
$$
$$
L_{0}u( x)=(1/2) a^{ij}( x) u_{x^{i}x^{j}} ( x)
$$
acting on functions given on $\bR^{d}$. 

For balls $B$ denote
$$
\osc (a,B )=
 |B |^{-2}
 \int_{y,z\in B }|a( y)-a( z)|\,dydz,
$$
$$
a^{\# }_{r}= \sup_{\substack {\rho\leq r\\ B\in\cB_{\rho}}}
\osc (a,B ) .
$$
 In the rest of the section we impose the following.
\begin{assumption}  
                                               \label{assumption 2.20.1} 
  
We have $p\in(1,\infty)$ and
 $a^{\#}_{R_{0}}\leq\theta_{0}=\theta_{0}(d,\delta,p)$, where $\theta_{0}>0$
is taken in a way to accommodate  Lemmas 3.3 and 3.4 of \cite{KrStr}.
\end{assumption}

\begin{remark}           \label{remark 6.30.1}
By Poincar\'e's inequality, for $B\in\cB_{\rho}$
$$
\osc (a,B )\leq N(d)\rho
\dashint_{B }|Da|\,dy
\leq N(d,\delta)\rho
\dashint_{B }|D\sigma|\,dy,
$$
so that   requiring in the future $\|  D\sigma\|$
to be sufficiently small will make $a$
satisfy Assumption 
\ref{assumption 2.20.1}.
In this sense Assumption 
\ref{assumption 2.20.1} is weaker than Assumption \ref{assumption 11.3.1}.
 
\end{remark}

As a consequence of Assumption
\ref{assumption 2.20.1}, derived in \cite{KrStr},
we have the following. Introduce
 $$
\Gamma=\{\Re \lambda\geq\lambda_{0}\}
\cup\{\varepsilon_{0}|\Im \lambda|\geq-\Re \lambda+\mu_{0}\},
$$
with $\varepsilon_{0}>0$ and $\mu_{0}>  \lambda_{0} $ which
depend only on $ d,\delta,p,R_{0}$
taken from Lemma 3.4 of \cite{KrStr}.
Then this lemma is the following.
\begin{lemma}
                         \label{lemma 11.5.2}
There exist
  $\mu_{0}>\lambda_{0}\geq 1, N_{0}$,
depending only on $d,\delta,p,R_{0}$, 
  such that, for any $u\in W^{2}_{p} $
and   $\lambda\in \Gamma$ we have
\begin{equation}
                                 \label{6.4.5}
 \| u_{xx}\|_{L_{p} } +|\lambda|\| u\|_{L_{p}  }
\leq N_{0} \|L_{0}u -\lambda u\|_{L_{p} }.
\end{equation}

\end{lemma}

To replace $L_{0}$ with $L$ we need to know
how to estimate the terms $b^{i}D_{i}u$.
Here is the Theorem from \cite{CF_90}
in which $d\geq 2$.
\begin{theorem}
                       \label{theorem 11.5.1}
Let $1<p<q\leq d$, $v$ be a function on $\bR^{d}$ such that for any $\rho>0$
and   $B\in\cB_{\rho}$  
$$
\Big(\dashint_{B}|v|^{q}\,dx\Big)^{1/q}
\leq \kappa\rho^{-1},
$$
where the constant $\kappa\in(0,\infty)$.
Then there exists a constant $N$,
depending only on $d,p,q$, such that
$$
\int_{\bR^{d}}|v(x)|^{p}|u(x)|^{p}\,dx
\leq N\kappa^{p} \int_{\bR^{d}}
| u_{x}(x)|^{p}\,dx
$$
as long as $u\in C^{\infty}_{0}$.
\end{theorem}

This theorem allows us to do the first step,
which generalizes
and implies Theorem \ref{theorem 11.5.1}
as $R_{0}\to\infty$.

\begin{lemma}   
                         \label{lemma 11.4.1}
For $1<p<q\leq d$ there exists a constant $N=N(p,q,  d )$ such that
\begin{equation}
                              \label{11.4.3}
\int_{\bR^{d}}|b(x)|^{p}|u(x)|^{p}\,dx
\leq N \| b \|^{p}
\int_{\bR^{d}} | u_{x}|^{p}\,dx
+N R_{0}^{-p }\|  b \|^{p} 
\int_{\bR^{d}} | u|^{p}\,dx
\end{equation}
as long as $u\in C^{\infty}_{0}$.
\end{lemma}

 Proof.  
Take $\zeta\in C^{\infty}_{0}(B_{R_{0}})$,
$\zeta\geq0$, such that
\begin{equation}
                         \label{11.4.5}
\int_{B_{R_{0}}}\zeta^{2p}\,dx=1, \quad
\zeta+ R_{0}|D\zeta|\leq N(d)R_{0}^{-d/(2p)}.
\end{equation}
 We claim that for any $\rho>0$ and  
$B\in \cB_{\rho}$   we have
\begin{equation}
                              \label{11.4.4}
\Big(\dashint_{B } |b \zeta|^{q}\,dx\Big)^{1/q}\leq NR_{0}^{-d/(2p)}\|  b \|  \rho^{-1}.
\end{equation}
Indeed, if $\rho\leq R_{0}$ it suffices
to use that $\zeta\leq NR_{0}^{-d/(2p)}$.
In case $\rho>R_{0}$, it suffices to use that
$$
\dashint_{B } |b \zeta|^{q}\,dx=
N\rho^{-d}\int_{B }|b \zeta |^{q}\,dx
\leq 
NR_{0}^{-qd/(2p)}\rho ^{-d}\int_{B_{R_{0}}}|b    |^{q}\,dx.
$$
$$
=NR_{0}^{-qd/(2p)}R_{0}^{d}\rho ^{-d}
\dashint_{B_{R_{0}}}|b    |^{q}\,dx
\leq NR_{0}^{-qd/(2p)}R_{0}^{d}\rho ^{-d}\|b\|^{q}
R_{0}^{-q}
$$
$$
=NR_{0}^{-qd/(2p)}(R_{0}/\rho)^{d-q}\rho^{-q}
\|b\|^{q}\leq NR_{0}^{-qd/(2p)}\|b\|^{q}\rho^{-q}.
$$
Now, in light of \eqref{11.4.4} by   Theorem
  \ref{theorem 11.5.1}
 $$
\int_{\bR^{d}}|b \zeta|^{p}|u|^{p}\,dx
\leq NR_{0}^{-d/2}\|  b \|^{p}
\int_{\bR^{d}}|Du|^{p}\,dx,\quad u\in C^{\infty}_{0}.
$$
We plug in here $\zeta(\cdot+y)$ and
$\zeta(\cdot+y)u$ in place of $\zeta$
and $u$, respectively. Then we get
$$
\int_{\bR^{d}}\zeta^{2p}(\cdot+y)|b |^{p}
|u|^{p}\,dx\leq
NR_{0}^{-d/2}\|  b \|^{p}
\int_{\bR^{d}}\zeta^{p}(\cdot+y)|Du|^{p}\,dx
$$
$$
+NR_{0}^{-d/2}\|  b \|^{p}
\int_{\bR^{d}}|D\zeta(\cdot+y)|^{p}| u|^{p}\,dx.
$$
After integrating through with respect to $y$
and using that by H\"older's inequality 
and \eqref{11.4.5}
$$
\int_{\bR^{d}}\zeta^{p}\,dy\leq
NR_{0}^{d/2},\quad
\int_{\bR^{d}}|D\zeta  |^{p}\,dy\leq
NR_{0}^{d/2-p},
$$
we come to \eqref{11.4.3}.
This proves the lemma.

On the basis of Lemmas \ref{lemma 11.5.2}
and \ref{lemma 11.4.1} we can repeat
what was done in \cite{KrStr} and obtain
the first part of the following result about the full operator $L$.

\begin{theorem}
                                                  \label{theorem 6.4.1}
Let $p\in(1,q)$. Then under   Assumptions  \ref{assumption 3.1.1}
and \ref{assumption 2.20.1} there   
exists $\bar N \geq 1$
depending only on $d,\delta,p$,
and $q$, such that,
if 
\begin{equation}
                          \label{11.7.1}
\bar N \|b\|\leq1,
\end{equation}
 then for any $u\in W^{2}_{p}$
and  $\lambda \in\Gamma$ ($\Gamma$ is introduced before Lemma \ref{lemma 11.5.2}),
\begin{equation}
                                 \label{6.4.6}
 \| u_{xx}\|_{L_{p} } +|\lambda|\| u\|_{L_{p}}
\leq N  \|L u -\lambda u\|_{L_{p}} 
\end{equation}
with $N$ depending only on $d,\delta,p,q,R_{0}$.
Furthermore, for any $\lambda \in\Gamma$ and $f\in L_{p}$
there is a unique  $u\in W^{2}_{p}$ such that
$\lambda u-Lu=f$.
\end{theorem}

The ``existence'' part of this theorem, as usual,
is proved by the method of continuity.

\begin{remark}
The use of \eqref{11.7.1} has very much in common with the ``form bounded'' condition
from \cite{Ki_20}:
 \begin{equation}
                            \label{11.20.2}
\int_{\bR^{d}}|b\phi|^{2}\,dx
\leq \delta \int_{\bR^{d}}|D\phi|^{2} \,dx
  + c_{\delta} \int_{\bR^{d}}|\phi|^{2} \,dx
\quad\forall \phi\in C^{\infty}_{0}.
\end{equation}
If you take here $\phi=\phi(x/\rho )$,
where $\phi\in C^{\infty}_{0}$ and $\phi=1$
on $B_{1}$ and $\phi=0$ outside $B_{2}$, then you get
$$
\dashint_{B_{\rho}}|b|^{2}\,dx\leq
N(d)\delta\rho^{-2}+N(d)c_{\delta}.
$$
This means that \eqref{11.20.1} is satisfied
with $q=2$. Condition \eqref{11.20.2}
is used in \cite{Ki_20}, among very many other things, to construct weak solutions of \eqref{6.15.2}
with constant $\sigma^{k}$'s.
We need a stronger condition \eqref{11.20.1}
with $q>d/2+1$ to prove the existence
of strong solutions.
\end{remark}

Below in this section we assume that 
\eqref{11.7.1} holds and denote by $R_{\lambda}f$ the solution $u$ from Theorem
\ref{theorem 6.4.1}. Next, exactly as
in \cite{KrStr} for complex
$t$  in the sector $S:=\{|\Im t|<\varepsilon_{0}\Re t\}$ 
with $\varepsilon_{0}$   introduced before Lemma \ref{lemma 11.5.2}
set
\begin{equation}
                                         \label{6.8.2}
 T_{t}=\frac{1}{2\pi i}\int_{\partial \Gamma}
e^{tz}R_{z}\,dz,
\end{equation}
where the integral is taken in a counter clockwise direction.
Below in this section 
$$
p\in(1,q).
$$ 

Here is Theorem 3.8 of \cite{KrStr}.
\begin{theorem}
                                                 \label{theorem 6.15.1}

Suppose that Assumptions  \ref{assumption 3.1.1}
and \ref{assumption 2.20.1} are satisfied and 
\eqref{11.7.1} holds. We assert the following.

(i) Formula \eqref{6.8.2} defines $ T_{t}$ in $S$ as
an analytic  semigroup of bounded operators
in $L_{p}$ with norms bounded 
by a constant, depending only on $\varepsilon,d,\delta,p$, $R_{0}$,  
as long as
$t\in\{|\Im t|\leq \varepsilon \Re t, |t|\leq(\varepsilon_{0}-
\varepsilon)^{-1}\}$ for any given
$\varepsilon<\varepsilon_{0}$;

(ii) The infinitesimal generator of this semigroup is
$L$ with domain
$W^{2}_{p}$;

(iii) For $g\in W^{2}_{p}$ the
function  $ T_{t}g(x)$ is a unique solution of the problem
$$
\partial_{t}u(t,x)=Lu(t,x),\quad t>0,\quad \lim_{t\downarrow 0}
\|u(t,\cdot)-g\|_{L_{p}}=0
$$
 in the class of $u$ such that
$u(t,\cdot)\in W^{2}_{p}$ and (strong $L_{p}$-derivative)
$\partial_{t}u (t,\cdot) \in L_{p}$
 for each $t>0$;

(iv) For any $T\in(0,\infty)$ there is a constant $N$,
depending only on $T,d$, $\delta,p,q,R_{0}$,    such that for each
$t\in(0,T]$ and
$f\in L_{p}$
\begin{equation}
                                         \label{6.15.1}
\| T_{t}f\|_{W^{2}_{p}}\leq \frac{N}{t}
\|f\|_{L_{p}},\quad 
\|D T_{t}f\|_{L_{p}}\leq \frac{N}{\sqrt{t}}
\|f\|_{L_{p}}.
\end{equation}

\end{theorem}

As in Remark 3.9 of \cite{KrStr} the properties of $ T_{t}$ listed in Theorem \ref{theorem 6.15.1} allow us to assert that, if 
$p>d/2 $ and $f\in W^{2}_{p}$, then $ T_{t}f $   has a modification that is  bounded and continuous in $x$, which we
still call $ T_{t}f$. Also   $ T_{t}f\to
 T_{s}f$ as $t\to s$ in $W^{2}_{p}$, and $ T_{t}f(x)\to  T_{s}f(x)$ uniformly on $\bR^{d}$.
Therefore, $ T_{t}f(x)$ is a bounded continuous function
on $[0,T]\times\bR^{d}$ for any $T\in(0,\infty)$.

Moreover,   for $0<t\leq T$, $f\in L_{p}$,
$q>p>d/2$,
and any $x\in\bR^{d}$
\begin{equation}
                                             \label{6.27.6} 
| T_{t}f(x)|\leq \frac{N}{t^{d/(2p)}}\|f\|_{L_{p}},
\end{equation}
where $N$ depends only on $T,d$, $\delta,p,q,R_{0}$.
 
We also need an approximation result
which, however, requires special 
approximation of $b$ and in this respect is
more restrictive 
than Theorem 3.10 of \cite{KrStr}.
Let $\zeta\in C^{\infty}_{0}(B_{1})$
be nonnegative spherically symmetric
with integral equal to one.
Define $\zeta_{n}(x)=n^{d}\zeta(nx)$ and
$b_{n}=b*\zeta_{n}$.

\begin{lemma}
                       \label{lemma 11.7.2}
For  any ball $B$
of radius $\rho\leq R_{0} $ we have
\begin{equation}
                           \label{11.7.3}
\Big(\dashint_{B}\sup_{n\geq 1/R_{0}}|b_{n}|^{q}\,dx\Big)^{1/q}
\leq N(d,q )\|b\|\rho^{-1},
\end{equation}
where $N(d,q )\geq1$.
\end{lemma}

Proof. If $B=B_{\rho}(x_{0})$, then on $B$
$$
\sup_{n\geq 1/R_{0}}|b_{n}|\leq N(d)
M(|b|I_{B_{2R_{0}}(x_{0})}),
$$
where $Mf$ is the Hardy-Littlewood maximal
function of $f$. By what is shown in the proof
of Theorem 1 of \cite{CF_88} the left-hand
side of \eqref{11.7.3} is dominated by
a constant, depending only on $d,q$, times $\rho^{-1}$
and times
the supremum over $r>0$ and  
  $B\in \cB_{r}$   of
\begin{equation}
                             \label{11.16.1}
I:=r\Big(\dashint_{B}|b|^{q}I_{B_{2R_{0}}(x_{0})}\,dx\Big)^{1/q}.
\end{equation}

If $r\leq R_{0}$ then $I$ is dominated by $\|b\|$. For $r\geq R_{0}$
we have
$$
I\leq N(d)r^{1-d/q}\Big(\int_{B_{2R_{0}}(x_{0})}|b|^{q}I_{B_{2R_{0}}(x_{0})}\,dx\Big)^{1/q}
$$
$$
\leq N(d)r^{1-d/q}\sup_{x\in\bR^{d}}
\Big(\int_{B_{ R_{0}}(x )}|b|^{q} \,dx\Big)^{1/q}
\leq N(d)r^{1-d/q}R_{0}^{d/q-1}\|b\|
\leq N(d)\|b\|.
$$
 This proves the lemma.

\begin{corollary}
                    \label{corollary 11.7.1}

If $u,u_{x}\in L_{p}$, then
$$
\lim_{n\to\infty}\int_{\bR^{d}}|b_{n}-b |^{p}|u |^{p}\,dx=0.
$$
\end{corollary}

Indeed, $|b_{n}-b |^{p}|u |^{p}\to 0$ (a.e.)
and 
$$
\sup_{n\geq 1/R_{0}}|b_{n}-b |^{p}|u |^{p}
\leq N\sup_{n\geq 1/R_{0}}|b_{n} |^{p}|u |^{p},
$$
where the last expression is integrable
in light of Lemmas \ref{lemma 11.7.2}
and \ref{lemma 11.4.1}.
\begin{theorem}
                                                    \label{theorem 6.23.1}

Let $q>p>d/2$ and let $a_{n} $, $n=1,2,...$, have the same meaning as $a $. Suppose that, for each $n$,  they satisfy
Assumptions \ref{assumption 3.1.1}
(with the same $\delta $) and \ref{assumption 2.20.1}  
(with the same $\theta_{0}$).
Suppose that \eqref{11.7.1} is satisfied with
$N(d,q )\|b\|$ in place of $\|b\|$,
where $N(d,q )$ is taken from \eqref{11.7.3}.
 Assume that $a_{n}\to a$ (a.e.) and take $b_{n}$ 
introduced before Lemma \ref{lemma 11.7.2}. Denote by $ T _{n,t}$
the semigroups constructed on the basis of \eqref{6.8.2}
when $R_{z}$ is replaced with $R _{n,z}$ that is the inverse
operator to $z-L_{n}$, where $L_{n}=(1/2)a^{ij}_{n}D_{ij}+b^{i}_{n}D_{i}$.
Then for any $t>0$ and $f\in L_{p}$ we have
$ T_{n,t}f\to T_{t}f$ in $W^{2}_{p}$ and, in particular,
uniformly on $\bR^{d}$ as $n\to\infty$.
\end{theorem}

The proof of this theorem is identical to the
proof of Theorem 3.10 of \cite{KrStr} up to the point where we show that
$$
\lim_{n\to\infty}\|\,|b^{n}-b|\,|(R_{z}f)_{x }|\|_{L_{p}}=0,
$$
which this time follows from Corollary
\ref{corollary 11.7.1}.  

\begin{remark}
                    \label{remark 11.14.1}
If $a$ and $b$ are smooth and $b$ is bounded,
then for any $f\in C^{\infty}_{0}$ there is a classical solution $u(t,x)$ of the problem $\partial_{t}u= 
Lu$, $u(0,\cdot)=f$. This solution and its derivatives decrease exponentially fast as $|x|
\to \infty$ and have all other qualities
listed in
Theorem \ref{theorem 6.15.1} (iii).
Therefore, $u(t,x)=T_{t}f(x)$. This shows,
in particular, that $T_{t}$ is independent of $p$. Owing to the maximum principle valid for $u$, we have
$$
0\leq T_{t}f\leq \sup f
$$
if an addition $f\geq0$. In light of \eqref{6.27.6}, this also holds for any
$f\in L_{p}$. The semigroup property of $T_{t}$
and \eqref{6.27.6}
imply further that for $t\geq 1$
$$
|T_{t}f|\leq T_{t}|f|=T_{t-1}T_{1}|f|
\leq\sup T_{1}|f|\leq N\|f\|_{L_{p}}.
$$
Thus, for $t>0$, $f\in L_{p}$,
$q>p>d/2$,
and any $x\in\bR^{d}$
\begin{equation}
                          \label{11.14.4}
| T_{t}f(x)|\leq \frac{N}{(t\wedge 1)^{d/(2p)}}\|f\|_{L_{p}},
\end{equation}
where $N$ depends only on $ d$, $\delta,p,q,R_{0}$.

These conclusions we obtained if
the coefficients are smooth and bounded.
By using the approximation Theorem
\ref{theorem 6.23.1} and mollifying $a$
we get the same conclusions in the general case
provided that \eqref{11.7.1} is satisfied with
$N(d,q,R_{0})\|b\|$ in place of $\|b\|$,
where $N(d,q,R_{0})$ is taken from \eqref{11.7.3}.
\end{remark}

\mysection{Stochastic equations with smooth coefficients}

Here, in addition to
Assumptions \ref{assumption 3.1.1} and \ref{assumption 11.3.1}, we suppose that $\sigma$ and $b$ are smooth
and bounded. First take $\zeta\in C^{\infty}_{0}(\bR^{d+1})$, $\zeta=\zeta(t,x),t\in\bR,x\in\bR^{d}$, such that $0\leq\zeta\leq1$,
$\zeta=1$ on $(-1,1)\times B_{1}$, and, for given $(t_{0},x_{0})\in\bR^{d+1}$, consider the equation
\begin{equation}
                          \label{11.10.1}
x_{t}=x_{0}+\int_{0}^{t}\sigma(x_{s})\,dw_{s}
+\int_{0}^{t}\zeta(t_{0}+s,x_{s})b(x_{s})\,ds.
\end{equation}

According to an obvious  modification
of Theorem 1.1 of \cite{Kr_21_1}  
with $p_{0}=p= d_{0}$, $ q_{0}=q=\infty$ 
(see, for instance,
the proof of  Theorem 1.2 in the same paper),
there is a constant $\hat b=\hat b(d,\delta)\in(0,\infty)$ such that if
\begin{equation}
                          \label{11.10.2}
\Big(\dashint_{B}|b|^{d_{0}}\,dx\Big)^{1/d_{0}}\leq \hat b \rho^{-1}
\end{equation}
for any $\rho\in(0,R_{0}]$ and   $B\in \cB_{\rho}$, then for any $R\in (0,R_{0}]$, $B\in \cB_{R}$, $x\in\bR^{d}$,
and Borel $ f  \geq0$
\begin{equation}
                                              \label{3.30.3}
 E_{x }\int_{0}^{\tau_{B} }
f(x_{t})\,dt\leq
\hat N 
R^{2-d/d_{0} }\|f \|_{L_{d_{0}}(B) },
\end{equation}
where $\hat N$ depends only on $d$ and $\delta $, $\tau_{B} $  is the first exit time
of the solution $x_{t}$   of \eqref{11.10.1}   from $B$ and we use the symbol $E_{x}$ to indicate that we are dealing with solutions
of stochastic equations started at $x$.

One can take $\zeta(t/n,x/n)$ in place of $\zeta$ in the above arguments and, since
any cylinder for large $n$ is absorbed in $(-n,n)\times B_{n}$ 
where the coefficients of \eqref{11.10.1}
coincide with the ones in \eqref{6.15.2},
we convince ourselves that
\eqref{3.30.3} holds for 
solutions of \eqref{6.15.2} once
\eqref{11.10.2} is satisfied (and the coefficients are smooth). By plugging in
$|b|$ in place of $f$ in \eqref{3.30.3}
and using that 
$$
\Big(\dashint_{B}|b|^{d_{0}}\,dx\Big)^{1/d_{0}}\leq \Big(\dashint_{B}|b|^{q}\,dx\Big)^{1/q},
$$
we get that
$$
\bar b_{R_{0}}:=\sup_{\substack{\rho\leq R_{0} \\ \,B\in\cB_{\rho}}}
\frac{1}{\rho}
\sup_{  x\in\bR^{d }}E_{ x}\int_{0}^{\tau_{B}}
|b( x_{s})|\,ds\leq N(d,\delta)\|b\|.
$$
It follows that there is a constant
$\bar N=\bar N(d,\delta)\geq1$ such that, if
\begin{equation}
                                     \label{12.18.3}
\bar N \|b\|\leq 1,
\end{equation}
then  \eqref{11.10.2} holds
for any $\rho\in(0,R_{0}]$ and   $B\in \cB_{\rho}$ and $\bar N\bar b_{R_{0}}
\leq 1$, where this $\bar N=\bar N(d,\delta)$
is taken from Theorem 2.3 of \cite{Kr_21_1}.
In this case  all conclusions of Theorem 2.3 of \cite{Kr_21_1}
hold true for any $R\leq R_{0}$.
Everywhere below in this section we suppose that \eqref{12.18.3} holds.
We thus have the following result, in which
\begin{equation}
                                                  \label{1.6.1}
\tau _{R} =\inf\{t\geq0: x_{t}\not\in B_{R} \},\quad
\gamma_{R} =\inf\{t\geq0:x_{t} \in \bar B_{R} \}.
\end{equation}
 
\begin{theorem}
                                      \label{theorem 8.2.1}
There is a
   constant $\bar \xi=\bar \xi(d,\delta)\in (0,1) $ 
such that   for any $R\leq R_{0}$ and $x$
\begin{equation}
                                          \label{8.2.2} 
  P_{ x}( \tau_{R}  \geq   R^{2} )\leq 1-\bar\xi.
\end{equation}
Moreover for $n=1,2,...$  
\begin{equation}
                                          \label{1.3.1} 
P_{ x}(\tau _{R} \geq nR^{2})
 \leq (1-\bar\xi)^{n},   
\end{equation}
so that 
\begin{equation}
                                          \label{3.7.2}
E_{ x}\tau _{R} \leq N(d,\delta)R^{2}.
\end{equation}
 
Furthermore,  the probability
starting from a point in  $\bar B_{9R/16}$
to reach the ball $\bar B_{R/16}$ before exiting
from $B_{R}$ is bigger than $\bar\xi$:
for any $x $ with $|x|\leq 9R/16$ 
\begin{equation}
                                          \label{1.2.1} 
P_{ x}(\tau _{R} >\gamma_{R/16} )\geq\bar\xi.   
\end{equation}

\end{theorem}

Once this result is established we can use all
results from \cite{KrPot} based on Theorem
2.3 from there. In particular, here is 
a particular case of Theorem 2.6 of \cite{KrPot}. We set $\tau'_{R}=\tau_{R}\wedge R^{2}$.

\begin{theorem}
                                     \label{theorem 8.20.1}
For any $\lambda,R>0$ satisfying $\lambda\geq
R_{0}^{-2}$ we have
\begin{equation}
                                          \label{8.20.1}
Ee^{-\lambda  \tau'_{R} }\leq
e^{\bar\xi/2}e^{- \sqrt{\lambda}  R
\bar\xi/2}.
\end{equation}

In particular, for  
  any   $R>0$ and $t\leq RR_{0}
\bar \xi/4 $ we have
\begin{equation}
                                             \label{10.2.2}
P(  \tau'_{R}  \leq  t )\leq 
 e^{\bar\xi/2}\exp\Big(-\frac{{\bar \xi}^{2}R^{2}}{16 t}\Big).
\end{equation}

\end{theorem}

This theorem implies Corollary 2.8 of \cite{KrPot}.
\begin{corollary}
                                        \label{corollary 10.26.1}
For any $m>0$ and
  $0\leq s\leq t $ we have
\begin{equation}
                                            \label{10.28.2}
E\sup_{r\in[s,t]}|x_{r}-x_{s}|^{ m}
\leq N(|t-s|^{ m/2}+|t-s|^{ m}),
\end{equation}
where $N=N(m,R_{0},\bar\xi)$.
\end{corollary}

 Since $\sigma$ and $b$ are smooth,
from the classical theory we know that
$E_{x}f(x_{t})= T_{t}f(x)$ for any
$f\in L_{p}$ with $p\in[1,\infty]$.
In particular,  
  \eqref{11.14.4} implies that
for $q>p>d/2$, $t>0$
\begin{equation}
                        \label{11.13.1}
E_{x}|f(x_{t})|\leq N (t\wedge 1)^{-d/(2p)}
\|f\|_{L_{p}},
\end{equation}
 where $N$ depends only on $d,\delta,p,q, R_{0}$.
As an obvious consequence of this estimate we also have that for any $q>p>d/2$, $\lambda\geq 1$
$$
E_{x}\int_{0}^{\infty}e^{-\lambda t}|f(t,x_{t})|
\,dt=\int_{0}^{\infty}e^{-\lambda t}
E_{x}|f(t,x_{t})|\,dt
$$
\begin{equation}
                        \label{11.13.2}
\leq N\int_{0}^{\infty}e^{-\lambda t}(t\wedge 1)^{-d/(2p) }\|f(t,\cdot)\|_{L_{p}}\,dt
 \leq N \lambda^{ (d+2)/(2p)-1}
\|f\|_{L_{p}(\bR^{d+1})},
\end{equation}
 where $N$ depends only on $d,\delta,p,q,R_{0}$.

\mysection{Properties of admissible solutions}

Recall that
$$
W^{1,2}_{p}([0,T]\times\bR^{d})=\{u:
u,u_{x},u_{xx},\partial_{t}u\in 
L_{p}([0,T]\times\bR^{d})\}.
$$
It is known that if $u\in W^{1,2}_{p}([0,T]\times\bR^{d})$
and $p>d/2+1$, then $u$ has a modification which
is bounded and continuous on $[0,T]\times\bR^{d}$.
Therefore, talking about $u$ of class $W^{1,2}_{p}([0,T]\times\bR^{d})$
we will always mean this modification. 
  Below by $x_{t}$ we mean
an admissible solution of \eqref{6.15.2},
corresponding to a $p\in(d/2+1,q)$,
starting at $x_{0}$ and
assuming that it exists.

\begin{theorem}[It\^o's formula]
                                             \label{theorem 6.7.1}
Let    
$u\in W^{1,2}_{p}([0,T]\times\bR^{d})$. Then  
with  probability one for all $t\in[0,T]$ we have
\begin{equation}
                                                     \label{6.16.2}
u(t,x_{t})=u(0,x_{0})+\int_{0}^{t}(\partial_{t}+L) u(s,x_{s})\,ds
+\int_{0}^{t}\sigma^{ik}D_{i}u(s,x_{s}) \,dw^{k}_{s},
\end{equation}
where the stochastic integral is a square integrable
martingale on $[0,T]$.
\end{theorem}

This theorem is proved by using \eqref{11.14.5}
in the same way as Theorem 1.3 of \cite{Kr_19_1}
is proved on the basis of Theorem 2.6 of \cite{Kr_19_1}. We only outline the main points
which are
\begin{equation}
                           \label{11.14.1}
E\int_{0}^{T}|b^{i}(x_{t})D_{i}u(t,x_{t})|\,dt
\leq N\|D^{2}u\|_{L_{p}((0,T)\times \bR^{d})},
\end{equation}
\begin{equation}
                           \label{11.14.10}
E\int_{0}^{T}| D u(t,x_{t})|^{2}\,dt
\leq N\|\partial_{t}u,D^{2}u\|_{L_{p}((0,T)\times \bR^{d})},
\end{equation}
where the constants $N$ are independent of $u$.

Estimate \eqref{11.14.1} immediately follows from \eqref{11.14.5} and Theorem \ref{theorem 11.5.1}. To prove \eqref{11.14.10} observe that
$$
1-(d+2)\Big(\frac{1}{p}-\frac{1}{2p}\Big)\geq 0
$$
so that
by embedding theorems (see, for instance,
Lemma 2.3.3 in \cite{LSU_68})
$$
\|\,|Du|^{2}\|_{L_{p}((0,T)\times \bR^{d})}
=\| Du \|_{L_{2p}((0,T)\times \bR^{d})}^{2}
\leq N\|\partial_{t}u, D^{2}u,u \|_{L_{ p}((0,T)\times \bR^{d})}^{2} .
$$
This and \eqref{11.14.5} imply \eqref{11.14.10}.

Here is a modification of Theorem 4.4 of
\cite{KrStr} in our situation.

\begin{theorem} 
                                             \label{theorem 6.16.2}
 Let $T\in(0,\infty)$  and 
$f\in L_{p} \cap L_{2p}$. Then 

(i) For each  $t>0$   we have 
$E  f(x_{t}) = T_{t}f(x_{0})$.
In particular,
\begin{equation}
                                         \label{12.2.1} 
   E |f(x_{t})|\leq N (t\wedge 1)^{-d/(2p)}
\|f\|_{L_{p}},
\end{equation}
where $N$ depends only on $d,\delta,p,q, R_{0}$;

(ii) For each  $t>0$ 
with  probability one   we
have
\begin{equation}
                                                     \label{6.16.3}
 f(x_{t})=   T_{t}f(x_{0})+\int_{0}^{t}\sigma^{ik}D_{i}
   T_{t-s}f(x_{s})\,dw^{k}_{s},
\end{equation}
where $\sigma^{ik}D_{i}
   T_{t-s}f(x )=\big(\sigma^{ik}D_{i}
   T_{t-s}f\big)(x )$ and similar notation is also used below;

(iii) For each  $t>0$   
\begin{equation}
                                                     \label{6.17.1}
   T_{t}f^{2}(x_{0})=(   T_{t}f(x_{0}))^{2}+
\sum_{k}\int_{0}^{t}   T_{s}\Big[\Big(\sum_{i}
\sigma^{ik}D_{i}   T_{t-s}f\Big)^{2}\Big](x_{0})\,ds.
\end{equation}

\end{theorem}

Proof. If $f\in W^{2}_{p}$, then $u(s,x):= T_{t-s}f(x)$, $s\leq t$,
satisfies the condition of Theorem \ref{theorem 6.7.1}
and we get \eqref{6.16.3} 
by that theorem.
By taking the expectations of both sides we get 
that $E f(x_{t})=     T_{t}f(x_{0})$.  
Then \eqref{12.2.1} follows from \eqref{11.14.4},
 By taking the expectations of
the squares of both sides of \eqref{6.16.3} we obtain
\eqref{6.17.1}. Thus, all assertions of the theorem
are true if $f\in W^{2}_{p}$.

Assertion (i) holds for any  $f\in L_{p}$,
which is seen from the fact that 
both $ T_{t}f(x_{0})$ and $ E f(x_{t})$
are bounded linear functionals on  a dense
subset  $ W^{2}_{p}$ of $L_{p}$.

Then, as $f^{n}\in  W^{2}_{p}$ tend to $f$
in $L_{p}\cap L_{2p}$, $   T_{t-s}f^{n}
\to    T_{t-s}f$ in $W^{2}_{p}$ for $s<t$ (see
\eqref{6.15.1}). By embedding theorems ($p\geq d/2$)
$D   T_{t-s}f^{n}
\to D   T_{t-s}f$ in $L_{2p}$ and in light of \eqref{11.14.4}
$$
   T_{s}\Big[\Big(\sum_{i}
\sigma^{ik}D_{i}   T_{t-s}f^{n}\Big)^{2}\Big](x_{0})
\to    T_{s}\Big[\Big(\sum_{i}
\sigma^{ik}D_{i}   T_{t-s}f\Big)^{2}\Big](x_{0})
$$
for any $0<s<t$. 
Furthermore, $(f^{n})^{2}\to f^{2}$ in $L_{p}$ and, due to \eqref{11.14.4},
$T_{t}(f^{n})^{2}(x_{0})\to T_{t}f^{2}(x_{0})$.
It follows by Fatou's lemma
(and \eqref{6.17.1}) that
\begin{equation}
                                                     \label{6.17.2}
   T_{t}f^{2}(x_{0})\geq (   T_{t}f(x_{0}))^{2}+
\sum_{k}\int_{0}^{t}   T_{s}\Big[\Big(\sum_{i}
\sigma^{ik}D_{i}   T_{t-s}f\Big)^{2}\Big](x_{0})\,ds.
\end{equation}

Hence, the right-hand side of \eqref{6.16.3}
is well defined. Furthermore,
$$
E \Big|\int_{0}^{t}\sigma^{ik}D_{i}
   T_{t-s}f(x_{s})\,dw^{k}_{s}-
\int_{0}^{t}\sigma^{ik}D_{i}
   T_{t-s}f^{n}(x_{s})\,dw^{k}_{s}\Big|^{2}
$$
$$
=\sum_{k}\int_{0}^{t}   T_{s}\Big[\Big(\sum_{i}
\sigma^{ik}D_{i}   T_{t-s}(f-f^{n})\Big)^{2}\Big](x_{0})\,ds
$$
$$
\leq    T_{t}(f-f^{n})^{2}(x_{0})- (   T_{t}(f-f^{n})(x_{0}))^{2}
 \leq T_{t}(f-f^{n})^{2}(x_{0})=E|f(x_{t})-f^{n}(x_{t})|^{2},
$$
where the   first  inequality is due to \eqref{6.17.2}. The last expression  
tends to zero  in light of \eqref{12.2.1}, which allows us to get \eqref{6.16.3}
by passing to the limit in its version with $f^{n}$
in place of $f$. After that \eqref{6.17.1} follows as above.
The theorem is proved.

Now we iterate \eqref{6.16.3}
and by repeating literally what
is done in \cite{KrStr} we come to the following
conclusions in which (as in \cite{KrStr})
$$
Q^{k}_{t}f(x)=\sigma^{ik}(x)D_{i}T_{t}f(x),
$$
  for $s_{1},...,s_{n}>0$  
we define
\begin{equation}
                                                                \label{6.26.1}
Q_{s_{n},...,s_{1}}f(x)=\sum_{k_{1},...,k_{n }}
\big[
Q^{k_{n }}_{s_{n }}\cdot...\cdot
Q^{k_{1}}_{s_{1}}f\big]^{2}(x ),
\end{equation} 
and $\cF^{w}_{t}$ is the completion of $\sigma(w_{s}:s\leq t)$.

\begin{theorem}
                                             \label{theorem 6.18.1}
Let $f\in L_{p}\cap L_{2p}$, $t>0$. Then
$$
E \big(f(x_{t})\mid\cF^{w}_{t}\big)=T_{t}f(x_{0})
$$
$$
+\sum_{m=1}^{\infty}\int_{t>t_{1}>...>t_{m}>0}T_{t_{m}}
Q^{k_{m}}_{t_{m-1}-t_{m}}\cdot...\cdot
Q^{k_{1}}_{t-t_{1}}f(x_{0})\,dw^{k_{m}}_{t_{m}}\cdot...\cdot dw^{k_{1}}_{t_{1}},
$$
where the series converges in the mean square sense.
\end{theorem}

\begin{theorem}
                                             \label{theorem 6.18.2}
Let $f\in L_{p}\cap L_{2p}$, $t_{0}>0$. Then
$f(x_{t_{0}})$ is $\cF^{w}_{t_{0}}$-measurable iff
\begin{equation}
                                                                \label{6.18.10}
\lim_{m\to\infty}  \int_{t_{0}>t_{1}>...>t_{m}>0}T_{t_{m }}
Q_{t_{m-1}-t_{m},...,t_{0}-t_{1}}f(x_{0})
\,d t_{m} \cdot...\cdot d t_{1}=0.
\end{equation}
Furthermore, under either of the above equivalent conditions
$$
 f(x_{t}) =T_{t}f(x_{0})
$$
\begin{equation}
                                                                \label{7.9.1}
+\sum_{m=1}^{\infty}\int_{t>t_{1}>...>t_{m}>0}T_{t_{m}}
Q^{k_{m}}_{t_{m-1}-t_{m}}\cdot...\cdot
Q^{k_{1}}_{t-t_{1}}f(x_{0})\,dw^{k_{m}}_{t_{m}}\cdot...\cdot dw^{k_{1}}_{t_{1}}.
\end{equation}

\end{theorem}

\begin{theorem}
                                             \label{theorem 6.18.3}
If equation \eqref{6.15.2} has two
admissible  solutions which
are not indistinguishable, then it does not have any admissible
strong solution. In particular, if \eqref{6.15.2} has an
  admissible strong solution, then it is
a unique admissible  solution.

\end{theorem}

\begin{theorem}
                                             \label{theorem 7.1.1}
If equation \eqref{6.15.2} has a strong
admissible
solution on one probability space then it has a strong admissible solution
on any other probability space carrying a $d_{1}$-dimensional
Wiener process.

\end{theorem}

Simple manipulations with \eqref{6.18.10}
as in \cite{KrStr} using \eqref{11.14.4}
lead to the following particular case of
Theorem 5.9 of \cite{KrStr}.

\begin{theorem}
                                           \label{theorem 6.26.1}
Let $f\in L_{p}\cap L_{2p}$. Then
$f(x_{t })$ is $\cF^{w}_{t }$-measurable for any $t>0$
if there exists a $\nu>0$ such that
\begin{equation}
                                                                \label{6.26.50}
  \Big\| \int_{\bR^{n }_{+} }e^{-\nu (s_{m-1}+...+s_{0})}
Q_{s_{m-1} ,...,s_{0} }f 
\,d s_{m-1 } \cdot...\cdot d s_{0}\Big\|^{p}_{L_{p}} \to 0
\end{equation}
as $m\to \infty$, where $\bR^{m}_{+}=(0,\infty)^{m}$.
\end{theorem}

We are going to prove that \eqref{6.26.50} holds under Assumptions \ref{assumption 11.3.1} and \ref{assumption 3.1.1} and
assuming that \eqref{11.15.1} holds
for an appropriate $N_{0}$,
by showing that the series composed of the left-hand sides
of \eqref{6.26.50} converges.

\mysection{Some estimates in the case of $C^{\infty}$ coefficients}
                                                \label{section 7.3.2}

We suppose that $\sigma^{k},b$ satisfy Assumption
\ref{assumption 3.1.1} and are infinitely differentiable
with each derivative bounded.

Let $(\Omega,\cF,P)$ be a complete probability space,
let $\{\cF_{t}\}$ be an increasing filtration of 
$\sigma$-fields $\cF_{t}\subset \cF$, that are complete.
Let 
$w_{t}$ be a $d_{1}$-dimensional Wiener process relative to
$\{\cF_{t}\}$. We also assume that there is a $(d+1)$ independent $d$-dimensional
Wiener, relative to  $\{\cF_{t}\}$, process $ B^{(0)}_{t},...,B^{(d)}_{t}$ 
independent
of $w_{t}$. Take $x,\eta\in \bR^{d}$, a nonnegative bounded
infinitely differentiable  $K_{0}$,   the role of which
will be emphasized later,  with each derivative bounded
given on $\bR^{d}$,
and consider the following system   
\begin{equation}
                                                        \label{6.20.3}
x_{t}=x+\int_{0}^{t}\sigma^{k}(x_{s})\,dw^{k}_{s}+
\int_{0}^{t}b(x_{s})\,ds,
\end{equation}
$$
\eta_{t}=\eta+\int_{0}^{t}\sigma^{k}_{(\eta_{s})}(x_{s})\,dw^{k}_{s}
+\int_{0}^{t}b_{(\eta_{s})}(x_{s})\,ds
$$
\begin{equation}
                                                        \label{6.20.4}
+\int_{0}^{t}K_{0}(x_{s})\,dB^{(0)}_{s}+\int_{0}^{t}K_{0}(x_{s}) \eta^{k}_{s} 
\,dB^{(k)}_{s} .
\end{equation}
As is well known, \eqref{6.20.3} 
 has a unique solution which we denote by $x_{t}(x)$.
By substituting it into \eqref{6.20.4} we see that the coefficients
of \eqref{6.20.4} grow linearly in $\eta$ and hence 
\eqref{6.20.4} also has a unique solution which we denote by
$\eta_{t}(x,\eta)$. By the way, observe that equation \eqref{6.20.4}
is linear with respect to $\eta_{t}$. Therefore
$\eta_{t}(x,\eta)$ is an affine function of $\eta$.
For the uniformity of notation we 
sometimes set $x_{t}(x,\eta)=x_{t}(x)$.  

For $t\geq0$ and $(x,\eta)\in\bR^{2d}$
consider the equation
$$
 \partial_{t}u(t,x,\eta)= (1/2)\sigma^{ik}\sigma^{jk}(x)u_{x^{i}x^{j}} (t,x,\eta)
+\sigma^{ik}\sigma_{(\eta)}^{jk}(x)u_{x^{i}\eta^{j}} (t,x,\eta)
$$
$$
+(1/2)\sigma_{(\eta)}^{ik} \sigma_{(\eta)}^{jk}(x)u_{\eta^{i}\eta^{j}}(t,x,\eta)
+(1/2)K^{2}_{0}(x)(1+|\eta|^{2})\delta^{ij}u_{\eta^{i}\eta^{j}}(t,x,\eta)
$$
\begin{equation}
                                                        \label{6.21.3}
+b^{i}(x)u_{x^{i}} (t,x,\eta)+b^{i}_{(\eta)}(x)u_{\eta^{i}} (t,x,\eta)
=:\check L(x,\eta)u(t,x,\eta)
\end{equation}
naturally related to system \eqref{6.20.3}-\eqref{6.20.4}.  
 
Here is Lemma 6.3 of \cite{KrStr}.

\begin{lemma}
                                                     \label{lemma 6.21.01}
 Let $x,\eta\in\bR^{d}$ and let $f (x)$
be infinitely differentiable with bounded derivatives.
 Then for any $t\in(0,\infty)$
$(t_{0}=t$)
$$
E\big[f_{(\eta_{t}(x,\eta))}(x_{t}(x))\big]^{2}
\geq\Big[(T_{t}f(x))_{(\eta)}\Big]^{2}
$$
\begin{equation}
                                                        \label{6.21.6}
+\sum_{m=1}^{\infty}\sum_{k_{1},...,k_{m}}
\int_{t>t_{1}>...>t_{m}>0}\Big[\big(T_{t_{m}}Q^{k_{m}}_{t_{m-1}-t_{m}}
\cdot...\cdot  Q^{k_{1}}_{t-t_{1}}f(x)\big)_{(\eta)}\Big]^{2}\,dt_{m}
\cdot...\cdot dt_{1}.
\end{equation}
\end{lemma}

Next, we want to estimate the left-hand side of \eqref{6.21.6}
which according to Lemma 6.1 of \cite{KrStr} satisfies
\eqref{6.21.3}. 

In the future we need a more precise
information than that provided
in Lemma 6.1 of \cite{KrStr}.   
\begin{lemma}
                      \label{lemma 11.16.1}
Take $f\in C^{\infty}_{0}$ and set
$$
u(t,x,\eta)=E\big[f_{(\eta_{t}(x,\eta))}(x_{t}(x))\big]^{2}.
$$
Then  $u$ is infinitely differentiable in $(x,\eta)$ and 
each of its derivatives is continuous in $t$ and 
$$
 |u(t,x,\eta)|+| u_{x}(t,x,\eta)|+
| u_{\eta}(t,x,\eta)|
$$
\begin{equation}
                                                        \label{11.16.60}
+| u_{xx}(t,x,\eta)|
+| u_{x\eta}(t,x,\eta)| +| u_{\eta\eta}(t,x,\eta)| 
\leq N e^{Nt-\kappa |x| } (1+|\eta|^{2}) ,
\end{equation}
where $N,\kappa>0$ are independent
of $x,\eta$.

\end{lemma}

Proof. We are going to use the terminology and results from Sections 2.7 and 2.8 of \cite{Kr_77}). Take unit $\mu,\nu\in \bR^{d}$.
As it follows from \cite{Kr_77}, the solution
$x_{t}(x)$ of \eqref{6.20.3} is infinitely
$LB$-differentiable in the direction of $\mu$
and the equations for the derivatives can be 
obtained by formal differentiation of \eqref{6.20.3}. This provides a sufficient
information to assert that the
solution $\eta_{t}(x,\eta)$ of \eqref{6.20.4}
is infinitely
$LB$-differentiable in the direction of $\mu$
in the variable $x$ and the equations for the derivatives can be 
obtained by formal differentiation of \eqref{6.20.4}. Similar assertion is true
for the derivatives of $\eta_{t}(x,\eta)$
with respect to $\eta$ in the direction of $\nu$ just because it
is   an affine function of $\eta$.
It follows, in particular, that
$u$ is infinitely differentiable in $(x,\eta)$.

By Theorem 2.8.8 of \cite{Kr_77} for any $T,r\in(0,\infty)$
\begin{equation}
                           \label{11.16.3}
E\Big(\sup_{t\leq T}\Big|LB-\frac
{\partial}{\partial \mu}x_{t}(x)\Big|^{r}
+\sup_{t\leq T}\Big|LB-\frac
{\partial^{2}}{\partial \mu^{2}}x_{t}(x)\Big|^{r}
\Big)\leq Ne^{NT},
\end{equation}
where $N$ is independent of $x,\eta,\mu,\nu$.
The derivative of $\eta_{t}(x,\eta)$ with respect to $\eta$ satisfies the same equation 
\eqref{6.20.4} but without the stochastic
integral of $K_{0}(x_{s})\,dB^{(0)}_{s}$.
Therefore this derivative admits an estimate
similar to \eqref{11.16.3}. Of course,
the second-order derivative of $\eta_{t}(x,\eta)$ with respect to $\eta$ is zero.
The mixed derivative
$$
LB-\frac
{\partial^{2}}{\partial \mu\partial \nu}\eta_{t}(x,\eta)
$$
satisfies the same equation as $\beta_{t}:=LB-(\partial/\partial \nu)\eta_{t}(x,\eta)$ but with zero initial data and a free term
$$
\int_{0}^{t}\sigma^{k}_{(\beta_{s})(\alpha_{s})}(x_{s})\,dw^{k}_{s}
+\int_{0}^{t}b_{(\beta_{s})(\alpha_{s})}(x_{s})\,ds
+\int_{0}^{t}K_{0(\alpha_{s})} (x_{s}) \beta^{k}_{s} 
\,dB^{(k)}_{s} ,
$$
where $\alpha_{t}:=LB-(\partial/\partial \mu)x_{t}(x )$. It follows very easily from 
\cite{Kr_77} that this derivative also admit
an estimate like \eqref{11.16.3}.

This and the fact that
$$
E\sup_{t\leq T}|\eta_{t,x}|^{r}
\leq N(1+|\eta|^{r})e^{NT}
$$
and $\sigma$ and $b$ are bounded
allows us to argue as  before
Theorem 6.4 of \cite{KrStr} and obtain
\eqref{11.16.60}
by using that $f$ has compact support. 
The lemma is proved.

In the future we might be interested in estimating not only
 the left-hand side of \eqref{6.21.6} but a slightly more
general quantity. Therefore, we take an infinitely
differentiable $f(x,\eta)\geq0$ such that
for an $m>0$ and a constant $N$
$$
\big(|f|+|f_{x}|+|f_{\eta}|+|f_{xx}|+|f_{x\eta}|+
|f_{\eta\eta}|\big)(x,\eta)\leq N(1+|\eta|)^{m}
$$
for all $x,\eta$  and such that $f(x,\eta)=0$
for all $\eta$ if $|x|\geq R$ for some $R>0$. Then denote 
$u(t,x,\eta)=  Ef[(x_{t},\eta_{t})(x,\eta)]$. According to \cite{KrStr}, there exist  constants  $ \mu >0$, $\kappa=\kappa(m)\geq0$, and a function
$M(t)$ bounded on each time interval $[0,T]$ such that for all $t $,
  $x,\eta$  we have
$$
 |u(t,x,\eta)|+| u_{x}(t,x,\eta)|+
| u_{\eta}(t,x,\eta)|
$$
\begin{equation}
                                                        \label{6.21.8}
+| u_{xx}(t,x,\eta)|
+| u_{x\eta}(t,x,\eta)| +| u_{\eta\eta}(t,x,\eta)| 
\leq M(t) e^{-\mu |x| } (1+|\eta|^{2})^{\kappa}.
\end{equation}
This justifies the integrations by parts we perform below.

Introduce
$$
h  = (1+| \eta|^{2} )^{-\kappa-d} 
$$ 
and observe that for a constant $N=N(d,\kappa)$ we have
$$
|\eta|\,|h _{\eta}|\leq N   h ,\quad 
 |((1+| \eta|^{2})
h)_{\eta \eta }|\leq N h  .
$$

\begin{theorem}
                                                   \label{theorem 6.21.1}
Let   $r\geq2$ and suppose that the above $u\geq0$.
 Then  there is a constant $\check N =
\check N 
(d,\delta,q,q_{0},\kappa) 
\geq1$  
such that if  
\begin{equation}
                                                        \label{7.1.1}
r\check N  (\|D\sigma\|+\|b\|)\leq 1,
\end{equation}
then
there exists  a constant  $N$,
depending only on    $d$,   $\delta$,    $q$, $q_{0}$, $\kappa$, $r$,
$R_{0}$,
  and there is a function $K_{0}$ 
such that for any $t\geq0$
\begin{equation}
                                                        \label{6.21.90}
\int_{\bR^{2d}}h ( \eta) u^{r}(t,x,\eta) \,dxd\eta\leq 
e^{Nt}\int_{\bR^{2d}}h ( \eta)f^{r}( x,\eta)\,dxd\eta.
\end{equation}

\end{theorem}

The proof of this theorem proceeds as usual
 by multiplying \eqref{6.21.3} by $h (\eta)u^{r-1}(t,x,\eta)$ 
and integrating by parts
over $[0,t]\times\bR^{2d}$. The integral of the left-hand side
is
$$
r^{-1}\int_{\bR^{2d}}h ( \eta)u^{r}(t,x,\eta)\,dxd\eta-r^{-1}
\int_{\bR^{2d}}h ( \eta)f^{r}( x,\eta)\,dxd\eta.
$$
Therefore, in light of Gronwall's inequality, to prove the theorem
it suffices to prove the following estimate.

\begin{lemma}
                                                        \label{lemma 6.21.5}
Let   $\kappa\geq0$, $r\in[2,\infty)$. 
Then  there is a constant 
$\check N \geq 1$  
depending only on    $d$,  $\delta$, $q$, $q_{0}$, $\kappa$,  
such that  if
\begin{equation}
                            \label{7.1.10}
r\check N (\|D\sigma\|+\|b\|)\leq 1,
\end{equation} 
then there exists a  constant  $N$,
depending only on    
$d$,   $\delta$,    $q$, $q_{0}$, $\kappa$, $r$,
$R_{0}$,
  and there is a function $K_{0}$  
  such that  for any   smooth 
function $v(x,\eta)\geq0$ (independent of $t$), for which condition
 \eqref{6.21.8} is satisfied with $v$ in place of $u$ and some $M$,
we have
\begin{equation}
                                                        \label{6.21.9}
\int_{\bR^{2d}}h ( \eta) v ^{r-1} ( x,\eta)\check Lv(x,\eta) \,dxd\eta
\leq N\int_{\bR^{2d}}h ( \eta) v^{r} ( x,\eta) \,dxd\eta .
\end{equation}

\end{lemma}

Proof. We basically repeat the proof
of Lemma 6.5 of \cite{KrStr} with
some changes caused by the weaker assumptions on $\sigma$ and $b$.
For simplicity of notation we drop the arguments $x,\eta$. 
We also write $U\sim V$ if
their integrals over $\bR^{2d}$ coincide, and $U\prec V$ if the integral of
$U$ is less than or equal to that of $V$.  Below the constants called $N$,
sometimes with indices, depend
 only on $d$,  $\delta$, $q$, $q_{0}$,  $\kappa$, $r$, $R_{0}$  unless specifically noted otherwise. Constants called $\hat N$
depend only on $d,\delta, q,q_{0}$, $\kappa$.

Set   $w=v^{r/2}$ and note simple formulas:
$$
v^{r-1}v_{x}=(2/r)ww_{x},\quad v^{r-2}v_{x^{i}}v_{x^{j}}
=(4/r^{2})w_{x^{i}}w_{x^{j}}.
$$
Then
denote by $\check L_{1}$ the sum of the first-order terms in $\check L$
 and observe that integrating by parts shows that
$$
h v^{r-1}  b^{i}_{(\eta)} v_{\eta^{i}}  \sim -(1/r)
h _{  \eta^{i}}  b^{i}_{(\eta)}v^{r}
-(1/r)h b^{i}_{x^{i}} v^{r}
$$
$$
\sim   (2/r)\eta^{k}h _{\eta^{i} } b^{i}w  w_{x^{k}} 
+(2/r)h b^{i}ww_{x^{i}} .
$$
Hence,
$$
hv^{r-1}\check L_{1}v \sim  
 (2/r)\eta^{k}h_{\eta^{i} } b^{i}w  w_{x^{k}} 
+(4/r)hb^{i}ww_{x^{i}} .
$$
Observe that  by Lemma \ref{lemma 11.4.1}
$$
\int_{\bR^{d}}|  b^{i}ww_{x^{k}}|\,dx \leq
\Big(\int_{\bR^{d}} |w_{x} |^{2}\,dx\Big)^{1/2}
\Big(\int_{\bR^{d}} |  b|^{2}|w|^{2}\,dx\Big)^{1/2}    
$$
\begin{equation}
                                                        \label{6.22.2}
\leq   \hat N \|b\| 
 \int_{\bR^{d}} |w_{x} |^{2}\,dx +
N \int_{\bR^{d}} |w  |^{2}\,dx ,
\end{equation}
where   $\hat N$ depends only on $d$, $q$
and $N$ depends only on $d$,   $q$, 
$R_{0}$, and, formally, $\|b\|$. 
But we suppress its dependence on $\|b\|$
because, in light of \eqref{7.1.10}
we assume from the start that $\|b\|,
\|D\sigma\|\leq 1$.
 
Since $|\eta|\,|h_{\eta}|\leq
N(\kappa,d)h$,
it follows that  
$$
\eta^{k}h_{\eta^{i} }  b^{i}ww_{x^{k}}
\prec \hat N  \|b\|   h|w_{x} |^{2}+Nh|w|^{2}.
$$ 
Similarly,
$(4/r)h  b^{i}ww_{x^{i}}\prec \hat N    \|b\|   h|w_{x} |^{2}+Nh|w|^{2}$ 
and we conclude that 
\begin{equation}
                                                        \label{6.22.1}
 h^{r}v^{r-1}\check L_{1}v\prec \hat N   \|b\|  h|w_{x} |^{2}+Nh|w|^{2}.
\end{equation}

Starting to deal with the second order derivatives note that  
$$
h v^{r-1}(1/2)\sigma^{ik}\sigma^{jk} v_{x^{i}x^{j}} \sim -
((r-1)/2)v^{r-2}h\sigma^{ik}v_{x^{i}} \sigma^{jk} v_{x^{j} } 
$$
$$
-(1/2)h\big[\sigma^{ik}_{x^{i}}\sigma^{jk}+
\sigma^{ik}\sigma^{jk}_{x^{i}}\big]v^{r-1}v_{x^{j}} =
-((2r-2)/r^{2}) h\sigma^{ik}w_{x^{i}} \sigma^{jk} w_{x^{j} } 
$$
$$
-(1/r)h\big[\sigma^{ik}_{x^{i}}\sigma^{jk}+
\sigma^{ik}\sigma^{jk}_{x^{i}}\big]ww_{x^{j}}  
\leq -(1/r) h\sigma^{ik}w_{x^{i}} \sigma^{jk} w_{x^{j} }
$$
$$
+h\Big|\big[\sigma^{ik}_{x^{i}}\sigma^{jk}+
\sigma^{ik}\sigma^{jk}_{x^{i}}\big]ww_{x^{j}}\Big|,
$$
where the inequality (to simplify the writing) is due to the fact that $r\geq2$.
In this inequality the first term on the right is dominated in the sense
of $\prec$ by 
$$
-(1/r)\delta h|w_{x}|^{2}
$$
 (see Assumption \ref{assumption 3.1.1}).
The remaining term contains $ww_{x^{i}}$ and we treat it as above. 
Then we get
\begin{equation}
                                                        \label{6.22.3}
h v^{r-1}(1/2)\sigma^{ik}\sigma^{jk} v_{x^{i}x^{j}}\prec
-\big[(1/r)\delta-\hat N \| D \sigma \| \big] h|w_{x}|^{2}+Nh|w|^{2}.
\end{equation}

Next,
$$
h v^{r-1}\sigma^{ik}\sigma_{(\eta)}^{jk} v_{x^{i}\eta^{j}}  \sim
-(r-1)h\sigma^{ik}v^{r-2} v_{\eta^{j}} \sigma_{(\eta)}^{jk} v_{x^{i}}
$$
$$
-v^{r-1} v_{x^{i}} \big[h_{\eta^{j}}  \sigma^{ik}\sigma_{(\eta)}^{jk}
+h \sigma^{ik}\sigma_{x^{j}}^{jk}]=-((4r-4)/r^{2})
h\sigma^{ik}  w_{\eta^{j}} \sigma_{(\eta)}^{jk} w_{x^{i}}
$$
$$
-(2/r)ww_{x^{i}} \big[h_{\eta^{j}}  \sigma^{ik}\sigma_{(\eta)}^{jk}
+h\sigma^{ik}\sigma_{x^{j}}^{jk}].
$$
We estimate the first term on the right roughly using
$$
|\sigma^{ik} w_{\eta^{j}} \sigma_{(\eta)}^{jk} w_{x^{i}}|\leq
\varepsilon|w_{x}|^{2}+\hat N\varepsilon^{-1}|\eta|\sum_{k}|\sigma^{k}_{x}|^{2}
|w_{\eta}|^{2}.
$$
The second term  contains $ww_{x^{i}}$ and allows the same handling as before.
 Therefore,
\begin{equation}
                                                        \label{6.22.4}
h v^{r-1}\sigma^{ik}\sigma_{(\eta)}^{jk} v_{x^{i}\eta^{j}} \prec
 (\varepsilon+ \hat N \|D\sigma \|)   h|w_{x}|^{2}+Nh|w|^{2}+\hat N\varepsilon^{-1}h|\eta|\sum_{k}|\sigma^{k}_{x}|^{2}
|w_{\eta}|^{2}.
\end{equation}

The last term in $h v^{r-1}\check  Lv$ containing $\sigma$ is
$$
h v^{r-1}(1/2)\sigma_{(\eta)}^{ik} \sigma_{(\eta)}^{jk} v_{\eta^{i}\eta^{j}} \sim
-((r-1)/2)h\sigma_{(\eta)}^{ik} v^{r-2} v_{ \eta^{j}}  \sigma_{(\eta)}^{jk}v_{\eta^{i}}
$$
$$
-(1/2)v^{r-1}\sigma_{(\eta)}^{ik} v_{\eta^{i} }\big[h_{\eta^{j}}  \sigma_{(\eta)}^{jk}
+h \sigma_{x^{j}}^{jk}\big]-(1/(2r))h (v^{r})_{\eta^{i}}
 \sigma_{x^{j}}^{ik} \sigma_{(\eta)}^{jk}
$$
$$
\prec \hat Nh(|\eta|^{2}|w_{\eta}|^{2}+w^{2})\sum_{k}|\sigma^{k}_{x}|^{2} +I,
$$
where  
$$
I=-(1/(2r))h(w^{2})_{\eta^{i}}
 \sigma_{x^{j}}^{ik} \sigma_{(\eta)}^{jk}
$$
$$
\sim
(1/(2r))w^{ 2}\sigma_{x^{j}}^{ik}\big[h_{\eta^{i}} \sigma_{(\eta)}^{jk}
+h\sigma_{x^{i}}^{jk}\big]\prec \hat Nh\sum_{k}|\sigma^{k}_{x}|^{2}w^{2}.
$$
To estimate the last term observe that by   Lemma \ref{lemma 11.4.1} 
\begin{equation}   
                                                          \label{7.5.2}
\int_{\bR^{d}}| \sigma^{k}_{x}|^{2}w^{2}\,dx\leq \hat N\|D\sigma\|^{2} \int_{\bR^{d}}|w_{x}|^{2}\,dx+N
\int_{\bR^{d}}|w |^{2}\,dx.
\end{equation} 
Above we had terms with $\|D\sigma\|$
and now we have   $\|D\sigma\|^{2}$.
To make formulas somewhat easier observe that
$\check N $, we are after, is bigger than one,
 so that $\|D\sigma\|\leq 1$ and
hence,
$$
I\prec \hat Nh\|D\sigma\|  |w_{x}|^{2}+Nh|w|^{2}
$$
and  
$$
h v^{q-1}(1/2)\sigma_{(\eta)}^{ik} \sigma_{(\eta)}^{jk} v_{\eta^{i}\eta^{j}}
\prec \hat Nh |\eta|^{2}|w_{\eta}|^{2} \sum_{k}|\sigma^{k}_{x}|^{2}
$$
\begin{equation}
                                                                 \label{7.9.3}
+ h w^{2}\Big(N+\hat N\sum_{k}|\sigma^{k}_{x}|^{2}\Big)+\hat Nh\|D\sigma\|  |w_{x}|^{2}.
\end{equation}

Finally,
$$
h v^{r-1}(1/2)K^{2}_{0}(1+|\eta|^{2})\delta^{ij}v_{\eta^{i}\eta^{j}} \sim
-((2r-2)/r^{2})hK^{2}_{0}(1+|\eta|^{2})|w_{\eta }|^{2}
$$
$$
-(2/r)K^{2}_{0}\big(h(1+|\eta|^{2})\big)
_{\eta^{i}}ww_{\eta^{i}}
$$
$$
\sim -((2r-2)/r^{2})hK^{2}_{0}(1+|\eta|^{2})|w_{\eta }|^{2}+(1/r) w^{2}
K_{0}^{2}\delta^{ij}\big(h(1+|\eta|^{2})\big)
_{\eta^{i}\eta^{j}}
$$
\begin{equation}
                                                                 \label{7.9.4}
\leq -(1/r)hK^{2}_{0}(1+|\eta|^{2})|w_{\eta }|^{2}+\hat N w^{2}
K_{0}^{2}h.
\end{equation}

By combining \eqref{6.22.1}, \eqref{6.22.3},
 \eqref{6.22.4}, \eqref{7.9.3}, and \eqref{7.9.4}, and using that
$|\eta|\leq 1+|\eta|^{2}$,
 we see that for any $\varepsilon\in(0,1]$

$$
h v^{q-1}\check Lv\prec  \Big[
\hat N_{1}   (\varepsilon+ \|b\|+\|D\sigma\|)
-\delta/r\Big]  h|w_{x} |^{2} 
$$
$$
 +\hat N_{2}\varepsilon^{-1}h(1+|\eta|^{2})\sum_{k}|\sigma^{k}_{x}|^{2}
|w_{\eta}|^{2} 
+N h w^{2} +\hat N_{3}h w^{2}\Big(K_{0}^{2}+\sum_{k}|\sigma^{k}_{x}|^{2}\Big) 
$$
\begin{equation}
                         \label{7.5.3}
-(1/r)hK^{2}_{0}(1+|\eta|^{2})|w_{\eta }|^{2}.
\end{equation}
 
  Here one sees clearly why introducing $K_{0}$,
which in no way helped us in \eqref{6.21.6}, is actually crucial.
With $K_{0}\equiv 0$ we would not be able
to estimate the term with $|w_{\eta}|^{2}$. 
Now, take and fix $\varepsilon$ so that $\hat N_{1} \varepsilon\leq \delta/(2r)$.
After that set
$$
K_{0}^{2}=1+\hat N_{2}r\varepsilon^{-1}\sum_{k}|\sigma^{k}_{x}|^{2}
$$
(1 is added to guarantee the smoothness of $K_{0}$)
and observe that according  to \eqref{7.5.2}
$$
\hat N_{3}h w^{2}
\Big(K_{0}^{2}+\sum_{k}|\sigma^{k}_{x}|^{2}\Big) =N  w^{2}
  h+\hat Nh w^{2}\sum_{k}|\sigma^{k}_{x}|^{2}
$$
$$
\prec \hat N_{4} h
\|D\sigma\|\,|w_{x}|^{2}+Nh w^{2}.
$$
Then \eqref{7.5.3} becomes
$$
h v^{r-1}\check Lv\prec  N  h|w|^{2}
-\Big[(1/(2r))\delta-(\hat N_{1}+\hat N_{4})    ( \|b\|+\|D\sigma\|) \Big] h|w_{x}|^{2}.
$$
We can certainly believe that $\hat N_{1}\geq1$, 
take $\check N$ in \eqref{7.1.10} to be equal to $(2 /\delta)(\hat N_{1}+\hat N_{4})$
($\geq1$),
and conclude that if \eqref{7.1.10}
holds, then
$$
h v^{r-1}\check Lv\prec  N  h|w|^{2}.
$$
The lemma is proved.

\mysection{Proof of Theorem \protect\ref{theorem 11.15.1}}

Set $p=(1/2)(d/2+1+q)$,
take $\zeta_{n}$ introduced before
Lemma \ref{lemma 11.7.2} and set $b_{n}=b*\zeta_{n}$, $\sigma_{n}=\sigma*\zeta_{n}$,
$a_{n}=\sigma_{n}\sigma^{*}_{n}$. Define
$$
\|D\sigma_{n}\|=\sup_{\substack{\rho\leq R_{0}
\\B\in \cB_{\rho}}}\rho\Big(\dashint_{B } |D\sigma_{n}  |^{q_{0}}dx\Big)^{1/q_{0}},
\quad \|b_{n}\|=\sup_{\substack{\rho\leq R_{0}
\\B\in \cB_{\rho}}}\rho\Big(\dashint_{B } |b_{n}|^{q }dx\Big)^{1/q }.
$$

\begin{lemma}
                        \label{lemma 11.15.4}
There is a constant $N_{0}=N_{0}(d,\delta,q_{0},q)\geq1$
such that if \eqref{11.15.1} is  
is satisfied with this $N_{0} $, then
for sufficiently large $n$

a) We have
\begin{equation}
                            \label{11.16.6}
p\check N (\|D\sigma_{n}\|+\|b_{n}\|)\leq 1,
\end{equation}
where $\check N=\check N(d,\delta/2,q,q_{0},2)$ 
is taken from Theorem \ref{theorem 6.21.1};

b) We have $a_{n,R_{0}}^{\#}\leq \theta_{0}
(d,\delta/2,p)$ and $\bar N(d,\delta/2,p,q)
N(d,q)\|b_{n}\|\leq 1$, where $\theta_{0}$
is taken from Assumption \ref{assumption 2.20.1}, $\bar N$ is the maximum of
$\bar N(d,\delta/2,p,q)$ from 
\eqref{11.7.1} and $\bar N(d,\delta/2)$
from \eqref{12.18.3}, and $N(d,q)$ is taken from
\eqref{11.7.3};   

c) The eigenvalues of $a_{n}$ are between
$\delta/2$ and $2\delta$.
\end{lemma}  

Proof. a) The possibility to find $N_{0}
=N_{0}(d,q,\check N)$ such that, \eqref{11.15.1}
would imply that
$q\check N  \|b_{n}\| \leq 1/2$, follows
from Lemma \ref{lemma 11.7.2}. This lemma
has an obvious counterpart applicable to
$D\sigma$ and this proves a).

b) The above argument and 
Remark \ref{remark 6.30.1} 
also take  care of b).

c) Denote by $\sigma$ the $d\times d_{1}$-matrix
whose columns are the $\sigma^{k}$'s and observe that
$$
|\sigma^{*}_{n}( x)\lambda|\leq \zeta_{ n}
( x) * |\sigma^{*}(x)\lambda|\leq\delta^{-1/2}|\lambda|.
$$
Therefore we need only prove that for sufficiently large $n$
\begin{equation}
                                                    \label{7.2.1}
|\sigma^{*}_{n}(x)\lambda|\geq |\lambda|\delta^{ 1/2}/\sqrt2.
\end{equation}
For any $y$ we have
$$
|\sigma^{*}_{n}(x)\lambda|\geq |\sigma^{*}(y)\lambda|-
 |(\sigma^{*}_{n}(x)-\sigma^{*}(y))\lambda|\geq |\lambda|\delta^{1/2}
-|(\sigma^{*}_{n}(x)-\sigma^{*}(y))\lambda| 
$$
$$
\geq |\lambda|\big(\delta^{1/2}
-| \sigma^{*}_{n}(x)-\sigma^{*}(y) |\big)
$$
Furthermore,
$$
\int_{\bR^{d}}| \sigma^{*}_{n}(x)-\sigma^{*}(x-y)  |\zeta_{n}(y)\,dy
$$
$$
\leq \int_{B_{1}}\int_{B_{1}}
| \sigma^{*}( x-z/n)-\sigma^{*}(x-y/n)  |\zeta (y)\zeta (z)\,dydz 
$$
$$
\leq N(d,q_{0})\|D\sigma\|,
$$
where the last inequality is due to
Poincar\'e.
 We see that to obtain c) it suffices to
have an appropriate $N_{0}=N_{0}(d,\delta,q_{0})$. The lemma is proved.

In the rest of the section we suppose that
\eqref{11.15.1} is satisfied with $N_{0}$
from Lemma \ref{lemma 11.15.4} and first prove
the existence of solutions.

\begin{theorem}
                     \label{theorem 11.13.1}
There exists a probability space and a $d_{1}$-dimensional
Wiener process on it such that 
equation \eqref{6.15.2}   has a solution
for which estimate \eqref{11.13.2} holds.
\end{theorem}

Proof. As usual we apply Skorokhod's
method. In light of Lemma \ref{lemma 11.15.4}, for sufficiently large $n$,
$\sigma_{n}$ and $b_{n}$ satisfy   
Assumptions \ref{assumption 3.1.1}, \ref{assumption 11.3.1} and \eqref{12.18.3} with $\delta/2$ in place of
$\delta$. Therefore, for the solutions
$x^{n}_{t}$ of 
\begin{equation}
                           \label{11.13.3}
x^{n}_{t}=x+\int_{0}^{t}\sigma_{n}(x^{n}_{s})
\,dw_{s}+\int_{0}^{t}b_{n}(x^{n}_{s})\,ds
\end{equation}
estimates \eqref{10.28.2} and \eqref{11.13.2} hold. After that we repeat the proof
of Theorem 2.6.1 of \cite{Kr_77} and see that
to finish proving the existence part of
the current theorem it suffices to show that
for any $T\in(0,\infty)$
\begin{equation}
                            \label{11.13.4}
 \int_{0}^{T}|b_{n}(x^{n}_{t})-b(x_{t})|
\,dt\to 0
\end{equation}
in probability as $n\to \infty$ provided that $x^{n}_{t}$
are solutions of \eqref{11.13.3} (with perhaps
different Wiener precesses for each $n$) and
$x_{t}$ is a continuous process such that
$x^{n}_{t}\to x_{t}$ in probability
for any $t\in[0,\infty)$.

Due to the convergence  of $x^{n}_{t}$
to $x_{t}$ estimate \eqref{11.13.2} holds
if $f$ is, in addition, bounded and continuous.
Then, of course, this estimate is extended
to all $f\in L_{p}$. Also obviously,
estimate \eqref{10.28.2} is true. This
shows that the probability of
$$
\{\sup_{t\leq T}|x^{n}_{t}|\geq R\}\cup
\{\sup_{t\leq T}|x _{t}|\geq R\}
$$
can be made as small as we like for all $n$
if $R$ is large enough. It follows that to prove
\eqref{11.13.4} it suffices to prove  that
\begin{equation}
                            \label{11.13.5}
\lim_{n\to\infty} E\int_{0}^{T}|\zeta(x^{n}_{t}) b_{n}(x^{n}_{t})-\zeta(x_{t}) b(x_{t})|
\,dt=0
\end{equation}
for any $\zeta\in C^{\infty}_{0}$.

Observe that for any bounded and continuous
$\bR^{d}$-valued $g$ the above limit
is dominated by 
$$
\nlimsup_{n\to\infty} E\int_{0}^{T}|\zeta(x^{n}_{t}) b_{n}(x^{n}_{t})-g(x^{n}_{t})|\,dt
+  E\int_{0}^{T}|g(x_{t})-\zeta(x_{t}) b(x_{t})|
\,dt,
$$
where both terms can be made as small as we like
because of estimate \eqref{11.13.2} valid for $x^{n}_{t}$ and $x_{t}$ and of the fact that $
\zeta b_{n}\to \zeta b$ in $L_{p}$
(even in $L_{q}$). This proves \eqref{11.13.4}
and establishes the existence of solution.
It turns out that in the above argument $x_{t}$
is exactly a solution for which, as we have seen, estimate \eqref{11.13.2} is valid.
The theorem is proved. 
 
Next, we prove that any admissible solution 
of \eqref{6.15.2} is strong.
Let $f\in C^{\infty}_{0}$.  First we deal with smooth
coefficients and develop necessary estimates.
Come back to Section \ref{section 7.3.2}
and consider the system \eqref{6.20.3}-\eqref{6.20.4} in which replace $\sigma,b$
with $\sigma_{n},b_{n}$ with $n$ so large that
the assertions a)-c) of Lemma \ref{lemma 11.15.4} are valid. Denote by $(x_{n,t},\eta_{n,t})
(x,\eta)$ the solution of the new system and
let $u_{n}(t,x,\eta)=E\big[f_{(\eta_{n,t}(x,\eta))}
(x_{n,t}(x))\big]^{2}$. Owing to 
Lemma \ref{lemma 11.16.1} estimate \eqref{6.21.8} holds with $\kappa=1$
and Lemma 
\ref{lemma 11.15.4} a) allows us
to use the conclusion of Theorem \ref{theorem 6.21.1} with $r=p$ and $u_{n}$ in place of $u$:
There exists $N=N(d,\delta,q,q_{0},R_{0})$
such that 
\begin{equation}
                                                        \label{11.17.1}
\int_{\bR^{2d}}h ( \eta) u^{p}_{n}(t,x,\eta) \,dxd\eta\leq 
e^{Nt}\int_{\bR^{2d}}h ( \eta)f^{p}( x,\eta)\,dxd\eta.
\end{equation}

 By Lemma \ref{lemma 6.21.01}  estimate
\eqref{11.17.1} implies that
 for $t\geq0$  
we have
\begin{equation}
                                                         \label{6.23.2}
\int_{\bR^{2d}}h ( \eta) v^{p}_{n}(t,x,\eta) \,dxd\eta\leq Ne^{Nt},
\end{equation}
where (and below) 
 $N$ depends only on $f$,  $d,\delta,q,q_{0}$, and $R_{0}$,
$$
v_{n} (t,x,\eta):=
$$
$$
\sum_{m=1}^{\infty}\sum_{k_{1},...,k_{m}}
\int_{t>t_{1}>...>t_{m}>0}\Big[\big(T_{n,t_{m}}Q^{k_{m}}_{n,t_{m-1}-t_{m}}
\cdot...\cdot  Q^{k_{1}}_{n,t-t_{1}}f(x)\big)_{(\eta)}\Big]^{2}\,dt_{m}
\cdot...\cdot dt_{1},
$$
and $T_{n,t},Q^{k}_{n,t}$ are constructed
from $\sigma_{n},b_{n}$ in the same way as $T_{t},Q^{k}_{t}$
are constructed from $\sigma,b$. By the  way this construction is possible thanks to 
Lemma \ref{lemma 11.15.4} b).

Obviously, $v_{n}(t,x,\eta)$ is a quadratic function
of $\eta$. Hence, \eqref{6.23.2} implies that, for any $R\in(0,\infty)$

\begin{equation}
                                                         \label{6.23.3}
\int_{\bR^{d}} \sup_{|\eta|\leq R} v_{n}^{p}(t,x,\eta) \,dx \leq Ne^{Nt}R^{2p}.
\end{equation}
Observe that in notation   \eqref{6.26.1}
naturally modified for $\sigma_{n},b_{n}$
$$
\sum_{k }v_{n}(t,x,\sigma^{k})
=\sum_{m=1}^{\infty} 
\int_{t>t_{1}>...>t_{m}>0}  Q_{n,t_{m},t_{m-1}-t_{n},
 ...,  t-t_{1}}f(x)   \,dt_{m}
\cdot...\cdot dt_{1}
$$
$$
=\sum_{m=1}^{\infty} 
\int_{S_{m}(t)}  Q_{n,s_{m},s_{m-1} ,  
 ...,s_{1},  t-(s_{1}+...+s_{m} )}f(x)   \,ds_{m}
\cdot...\cdot ds_{1}=:\sum_{m=1}^{\infty}I_{n,m}(t,x),
$$
where 
 $S_{m}(t)=\{(s_{1},...,s_{m}):s_{k}>0,s_{1}+...+s_{m}< t \}$.
Next, for $\nu>0$ by H\"older's inequality
$$
\sum_{m=1}^{\infty}\int_{\bR^{d}}\Big(\int_{0}^{\infty}
e^{-\nu t}I_{n,m}(t,x)\,dt\Big)^{p}\,dx
$$
$$
\leq \nu^{1-p}\int_{0}^{\infty}e^{-\nu t}\Big(\sum_{m=1}^{\infty}
\int_{\bR^{d}}I^{p}_{n,m}(t,x)\,dx\Big)dt
$$
$$
\leq \nu^{1-p}\int_{0}^{\infty}e^{-\nu t}\int_{\bR^{d}}
\Big(\sum_{k }v_{n}(t,x,\sigma^{k})\Big)^{p}\,dxdt ,
$$
which thanks to \eqref{6.23.3} implies that  for appropriate $\nu$,
depending only on $f$,  $d,\delta,q,q_{0}$, and $R_{0}$,
\begin{equation}
                                                              \label{6.27.1}
\sum_{m=1}^{\infty}\int_{\bR^{d}}\Big(\int_{0}^{\infty}
e^{-\nu t}I_{n,m}(t,x)\,dt\Big)^{p}\,dx
 \leq N,
\end{equation}
where   $N$ depends    only on $f$,  $d,\delta,q,q_{0}$, and $R_{0}$.

Now we let $n\to\infty$ in \eqref{6.27.1}.
Observe that since $\sigma_{n}\to\sigma$,
$b_{n}\to b$ (a.e.) we have $a_{ R_{0}}^{\#}\leq \theta_{0}
(d,\delta/2,p)$ and $\bar N(d,\delta/2,p,q)N(d,q)\|b \|\leq 1$, where $\theta_{0}$
is taken from Assumption \ref{assumption 2.20.1}, $\bar N$ is taken from 
\eqref{11.7.1}, and $N(d,q)$ is taken from
\eqref{11.7.3}. Therefore, the semigroup $T_{t}$ is well defined as in Section \ref{section 3.11.2}.

Also note that in light of Theorem \ref{theorem 6.23.1} for any $\eta\in\bR^{d}$,
$t>t_{1}....>t_{m}>0$
$$
\big(T_{n,t_{m}}Q^{k_{m}}_{n,t_{m-1}-t_{m}}
\cdot...\cdot  Q^{k_{1}}_{n,t-t_{1}}f(x)\big)_{(\eta)}\to
\big(T_{ t_{m}}Q^{k_{m}}_{ t_{m-1}-t_{m}}
\cdot...\cdot  Q^{k_{1}}_{ t-t_{1}}f(x)\big)_{(\eta)}
$$
in $L_{p}$. It follows by Fatou's lemma that
$$
\nliminf_{n\to \infty}I_{n,m}\geq I_{m},
\quad
\sum_{m=1}^{\infty}\int_{\bR^{d}}\Big(\int_{0}^{\infty}
e^{-\nu t}I_{ m}(t,x)\,dt\Big)^{p}\,dx<\infty.
$$

Finally, by observing that
$$
\int_{0}^{\infty}
e^{-\nu t}I_{m}(t,x)\,dt=\int_{\bR^{m+1}_{+}}e^{-\nu(s_{0}+...+s_{m})}
Q_{s_{m},...,s_{0}}f(x)\,ds_{m}\cdot...\cdot ds_{0}
$$
and referring to Theorem \ref{theorem 6.26.1},   we 
conclude that $f(x_{t})$ is $\cF^{w}_{t}$-measurable for any $t\geq0$.
The arbitrariness of $f$ and $t$
finishes the proof.

\end{document}